\documentclass[12pt,a4paper]{amsart}
\usepackage[utf8]{inputenc}	
\usepackage[T1]{fontenc}
\usepackage[in,headings]{fullpage}
\usepackage{amsbsy, amscd, amsfonts, amsmath, amsrefs, amssymb, amsthm}
\usepackage{graphicx}
\usepackage{float}
\usepackage{color}   
\usepackage[colorlinks=true,linkcolor=blue,citecolor=blue,urlcolor=blue]{hyperref}

\newtheorem{theorem}{Theorem}

\newtheorem{lemma}[theorem]{Lemma}

\theoremstyle{definition}

\theoremstyle{remark}

\newcommand{\R}{\mathbf{R}}

\renewcommand{\Re}{\mathop{\mathrm{Re}}\nolimits}
\renewcommand{\Im}{\mathop{\mathrm{Im}}\nolimits}
\newcommand{\Rzeta}{\mathop{\mathcal R }\nolimits}

\newfont{\cmbsy}{cmbsy10}
\newfont{\cmmib}{cmmib10}
\newcommand{\Orden}{\mathop{\hbox{\cmbsy O}}\nolimits}
\newcommand{\orden}{\mathop{\hbox{\cmmib o}}\nolimits}

\begin{document}

\title[Statistic of zeros of $\Rzeta$]
{Statistic of zeros of Riemann auxiliary function.}
\author[Arias de Reyna]{J. Arias de Reyna}
\address{%
Universidad de Sevilla \\ 
Facultad de Matem\'aticas \\ 
c/Tarfia, sn \\ 
41012-Sevilla \\ 
Spain.} 




\email{arias@us.es, ariasdereyna1947@gmail.com}


\begin{abstract}
We have computed all zeros $\beta+i\gamma$ of $\Rzeta(s)$ with $0<\gamma<215946.3$. A total of 162215 zeros with 25 correct decimal digits. 
In this paper we offer some statistic based on this set of zeros. Perhaps the main  interesting result is that $63.9\%$ of these zeros satisfies $\beta<1/2$. 
\end{abstract}

\maketitle

\section{Introduction} 

Riemann asserts in \cite{R} that $\zeta(s)$ have a number of zeros in the critical line to height $T$ similar to the total number of zeros to this height. In a letter to Weiertrass \cite{R2}*{p.~823} he repeats this and says that this is difficult to prove and depends on a new representation of the function $\Xi(t)$ he has not communicated. Siegel \cite{Siegel} tried to clear this, but concluded that in the remaining papers of Riemann there is no sign of proof of this result. Nevertheless he finds in Riemann's paper a function 
\[\Rzeta(s)=\int_{0\swarrow1}\frac{x^{-s}e^{\pi i x^2}}{e^{\pi i x}-e^{-\pi i x}}\,dx\]
which zeros are related to the zeros of $\zeta(s)$. Applying it, Siegel shows that $\zeta(s)$ have $>cT$ zeros on the critical line with $0<\gamma\le T$. 

We have studied this function and its zeros \cite{A86}, \cite{A166},  \cite{A98}, \cite{A100}, \cite{A108}. We have computed the zeros $\beta+i\gamma$ of $\Rzeta(s)$  with $0<\gamma<215946.3$. The results are contained in a file \texttt{Zeros25digits.txt}\footnote{Available with the tex source of this document.}, which contains 162215 zeros. The real and imaginary parts of each zero are given with 25 correct decimal digits. A glimpse of all these zeros can be seen in Figure \ref{F:plotzeros2}. In this paper, we give some statistics related to these zeros. 

In Section \ref{S:zeros} we give some graphic idea of the zeros and resume the results about its situation that we have proved.
The zeros $\rho=\beta+i\gamma$ of $\Rzeta(s)$ with $\gamma>0$  are distributed 
along a region near the imaginary axis. In Section \ref{S:order} I explain how I  presumably get all zeros to a certain height and how we ordered the zeros. The connection of the zeros of $\Rzeta(s)$ and those of $\zeta(s)$ is explained in Section \ref{S:zerosR-and-zeta}.

Siegel proves that the number of zeros of $\Rzeta(s)$ with $0<\gamma\le T$ is just half the number of zeros of $\zeta(s)$ with error $\orden(T)$. We explain in \ref{S:5} how I predicted an additional term in this approximation. Later we gave a proof of this new term in \cite{A185}.

Siegel asserted without proof that for any $\varepsilon >0$  the zeros with $\gamma>0$ satisfies $-C\gamma^\varepsilon\le\beta\le2$ for some constant. I have not been able to prove this. Only an inequality of type $-A\gamma^{2/5}\log\gamma\le \beta\le 2$ (see \cite{A98}). To study the question of the left limit of zeros, in Section \ref{S:6}, we consider the sequence of zeros 
$\rho_{n_k}$ such that $-\beta_{n_k}$ are records. It happens that the corresponding $\gamma_{n_k}\approx 2\pi (k+1)^2$. And this is the first appearance of a cyclic behaviour of the zeros of $\Rzeta(s)$. But the more interesting question of the $\beta_{n_k}$ is more elusive. We need more data to answer this question. 

One of the results of Siegel on \cite{Siegel}
is that 
$h(t)=-\sum_{0<\gamma\le T}\beta=\frac{T}{4\pi}\log 2+\orden(T)$.
In Sections \ref{S:7} and \ref{S:8} we conjecture an extra periodic term in this equation \eqref{E:conj}. 

\begin{figure}[H]
\begin{center}
\includegraphics[width=\hsize]{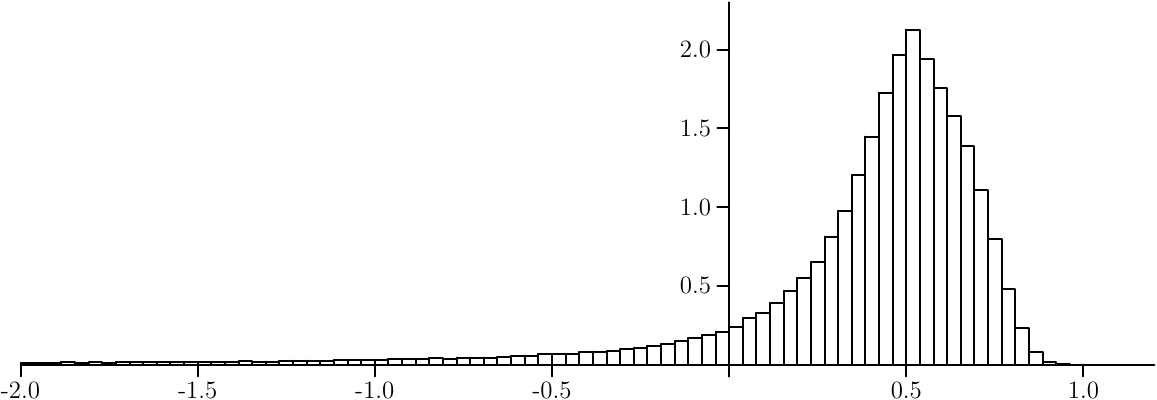}
\caption{Distribution horizontal of zeros}
\label{F: horizontal}
\end{center}
\end{figure}

Finally, in the last section we consider the horizontal distribution of the zeros. A histogram of the computed zeros is shown in Figure \ref{F: horizontal}.
In an appendix, we give some tables of data on this distribution. 
More or less 2/3 of the computed zeros satisfies $\beta<1/2$. If this can be proved, we will get $2/3$ of the zeros of $\zeta(s)$ on the critical line. This will improve on the results obtained until now by the refinements of Levinson's method. By  Theorems analogous to the density theorems \cite{A174} most of the zeros with $\beta>1/2$ are very close to the critical line. This makes it plausible that each of these zeros also contributes to two zeros of $\zeta(s)$ on the critical line.

\section{Situation of the zeros}\label{S:zeros}
Siegel \cite{Siegel}*{p.~299} wondered  where  the zeros of $\Rzeta(s)$ are situated.
His results in the paper,  some refinements, and new results obtained by us
give us some answers to this question. 

\begin{figure}[H]
\begin{center}
\includegraphics[width=0.7\hsize]{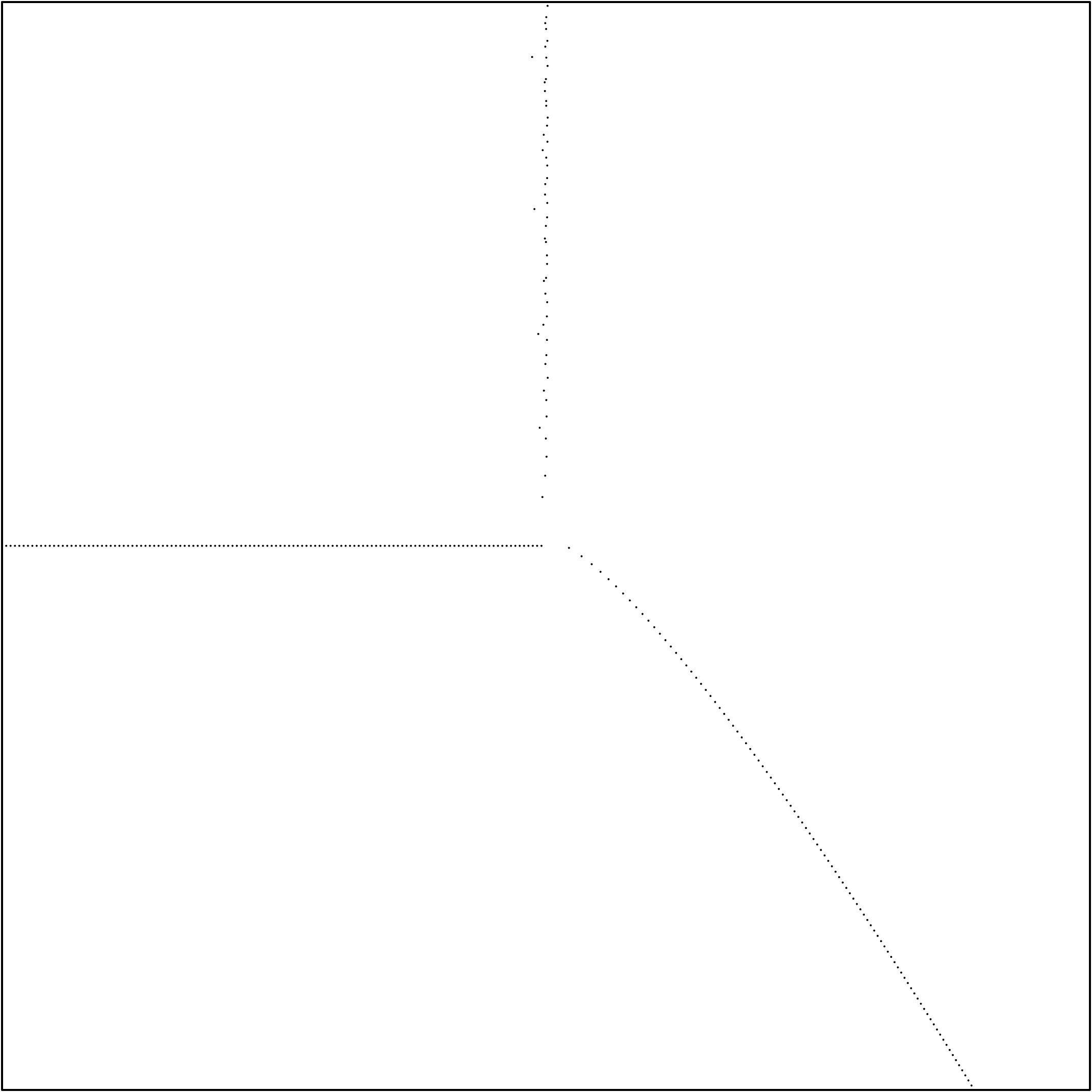}
\caption{Zeros of $\Rzeta(s)$ on the square $(-250,250)^2$}
\label{figzeros}
\end{center}
\end{figure}

We see in the figure that there are three lines of zeros. The trivial zeros at $-2n$ form one of these lines, then there are other ones appearing to be near the imaginary axis and finally there is a line of zeros on the fourth quadrant. 

\begin{figure}[H]
\begin{center}
\includegraphics[width=0.95\hsize]{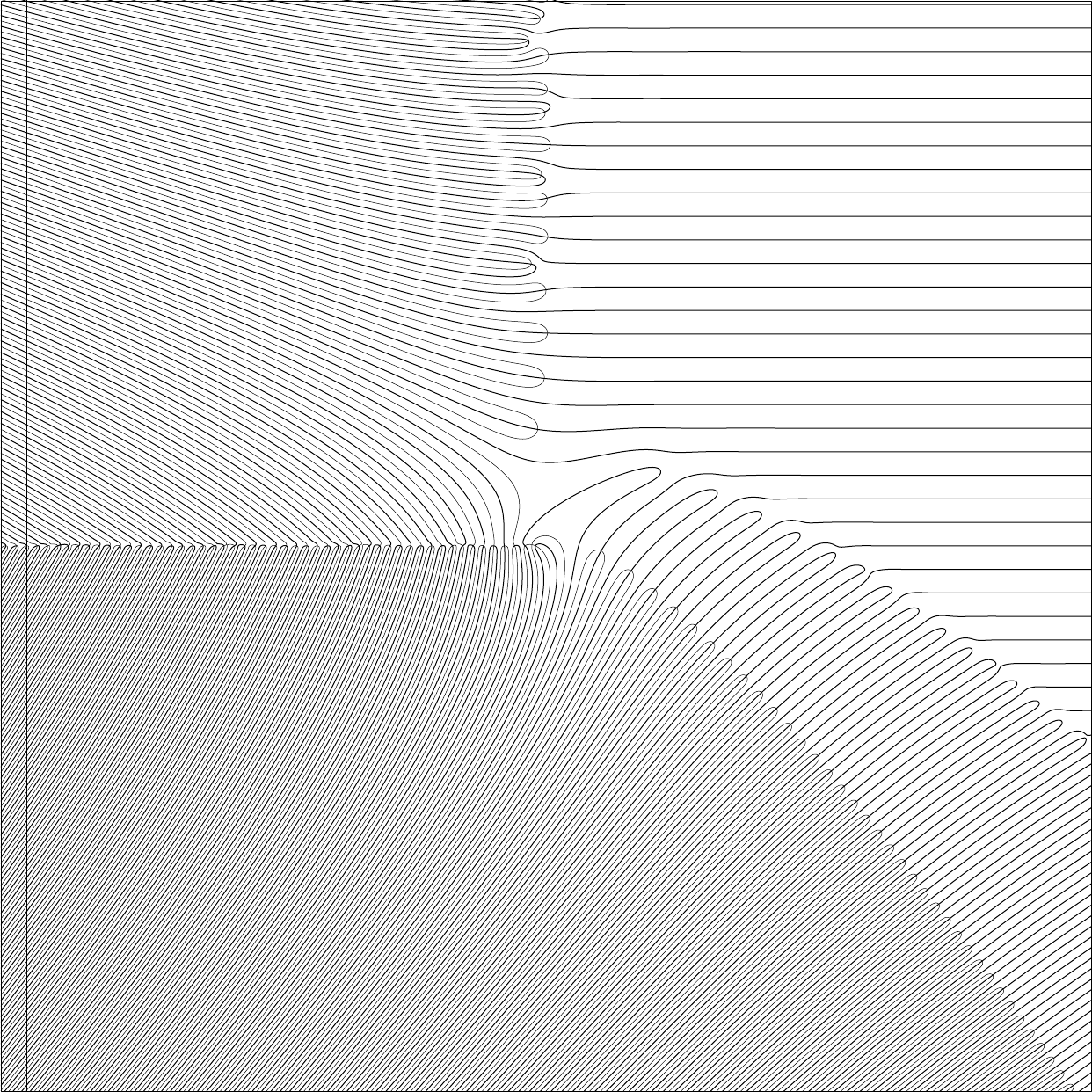}
\caption{x-ray of $\Rzeta(s)$ in $(-105,105)^2$ with the line of seeds marked.}
\label{xrayRzeta}
\end{center}
\end{figure}

In \cite{A100}*{Cor.~14} we show that there is no zero in the region limited by $\sigma\ge2$ and a line approximately parallel to the line of zeros in the fourth quadrant.  Theorem 16 in \cite{A100} proves that there are no zeros for a region that occupies most of the third quadrant and part of the fourth quadrant below the line of zeros on the fourth quadrant.  

Theorem 7 in \cite{A98} proves that the only zeros in a region to the left of the imaginary axis, but separated of it, are the trivial ones, which are simple.

The left limit of the upper zeros is difficult. Siegel asserts that for any $\varepsilon>0$ the zeros $\rho=\beta+i\gamma$ are to the right of the line $(-t^\varepsilon, t)$ so that $\beta+\gamma^\varepsilon>0$. But he only proves it for 
$\varepsilon=3/7$. In \cite{A98} we show that his reasoning only extends $\beta+A \gamma^{2/5}\log \gamma$. 

All our computed zeros with $\gamma>0$ satisfies $\beta<1$. We have proved \cite{A173} this is true for $\gamma\ge t_0$ for a high value of $t_0$, and we conjecture that this is true in general. 

All these results do not apply for zeros near the origin, contained in a circle of radius 10000 or so; see each particular result for details. This is because all are obtained from asymptotic expansions.
The computed zeros make clear that the results are also true for smaller  values. 

In \cite{A186} we show that all trivial zeros are simple. 

In \cite{A108} we give information about the zeros in the fourth quadrant. They are very regular and can be computed rather easily. I will concentrate here on the zeros with $\gamma>0$. A list of the first 2122 zeros with $\gamma<0$ can be found at the end of the TeX-file for \cite{A108}. 

\section{Connection with zeros of zeta}\label{S:zerosR-and-zeta}
Siegel showed  (see \cite{A166}*{eq.~(3.8)}, and Titchmarsh \cite{T}*{eq.~(4.17.2)})
\[Z(t)=2\Re\bigl\{e^{i\vartheta(t)}\Rzeta(\tfrac12+it)\bigr\}, \qquad t\in\R,\]
where, as usual, $\zeta(\frac12+it)=Z(t)e^{-i\vartheta(t)}$ with real analytic functions $Z(t)$ and $\vartheta(t)$.

This can also be written as 
\[-\frac{\Xi(t)}{\frac14+t^2}=\Re\bigl\{\pi^{-s/2}\Gamma(s/2)\Rzeta(s)\bigr\},\qquad s=\tfrac12+it.\]
Therefore, the zeros of $\zeta(s)$ on the critical line coincide with the cuts of the imaginary lines of the x-ray of the function $\pi^{-s/2}\Gamma(s/2)\Rzeta(s)$ with the critical line. 

\section{Order of the zeros}\label{S:order}

\begin{figure}[H]
\begin{center}
\includegraphics[width=0.6\hsize]{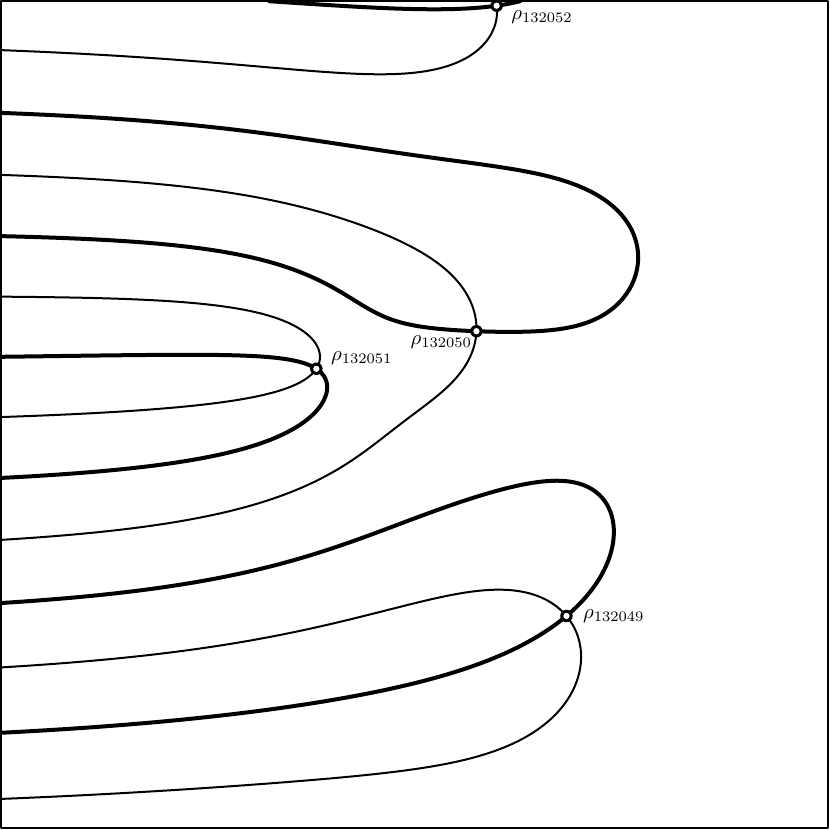}
\caption{Zero $\rho_{132050}$ is higher than $\rho_{132051}$.}
\label{orderzeros}
\end{center}
\end{figure}

We know that there are two constants $A$ and $t_0$ such that all zeros $\rho=\beta+i\gamma$ with $t_0<\gamma<T$ satisfies $\beta>1-At^{2/5}\log t$. We know that such a number exists.  The number of zeros with $0<\gamma\le t_0$ is finite. Therefore, there is some constant $\sigma_0$ such that all zeros with $0<\gamma<215946.3$ satisfies $\beta>\sigma_0$. Looking at the x-ray of $\Rzeta(s)$,
we see that the lower  zeros are each in a line that cuts the line $\sigma=-100$. Since the extreme zeros change gradually when we increase the height, when we arrive at the height $T=215946.3$ we can trust that each zero below this line satisfies $\beta>-100$.  Our computation of zeros starts from a seed, a point where $\Rzeta(-100+i t)$ is purely imaginary and such that 
$\Rzeta(-100+i t)/i<0$. We make a list of these seeds $0<t_1<t_2<\cdots$. (This imaginary line will cut $\sigma=-100$ at the other points where $\Rzeta(-100+i t)/i>0$). Each zero corresponds to a seed, and we find it following the line until we get $\Im\Rzeta(s)$ small. That is when we arrive at the proximity of the zero. Then we use the Newton method to get the zero with the required precision. 

The first seed is $t_1\approx60.969$. The first values of $t$ with $\Rzeta(-100+i t)/i<0$ are related to trivial zeros. %
There is a natural order between the seeds. And we will call zero $\rho_n=\beta_n+i\gamma_n$ the one related to the seed $t_n$. This gives us an order to the zeros that approximately coincides, but is not equal, to the order determined by the $\gamma_n$. 
Figure \ref{orderzeros} shows a case in which a higher zero has a lower number than the other zero.  In situations like this, it is common that if we follow the line from the seed $t_{132051}$ with not much precision (and this is very important to get a fast program), then the Newton method gives in the two cases the zero $\rho_{132050}$. We then have to  refine  the search for zeros in these cases later. It is not easy to get a complete list of zeros. And it is important to know which seed is failing in each case.

\begin{figure}[H]
\begin{center}
\includegraphics[width=0.9\hsize]{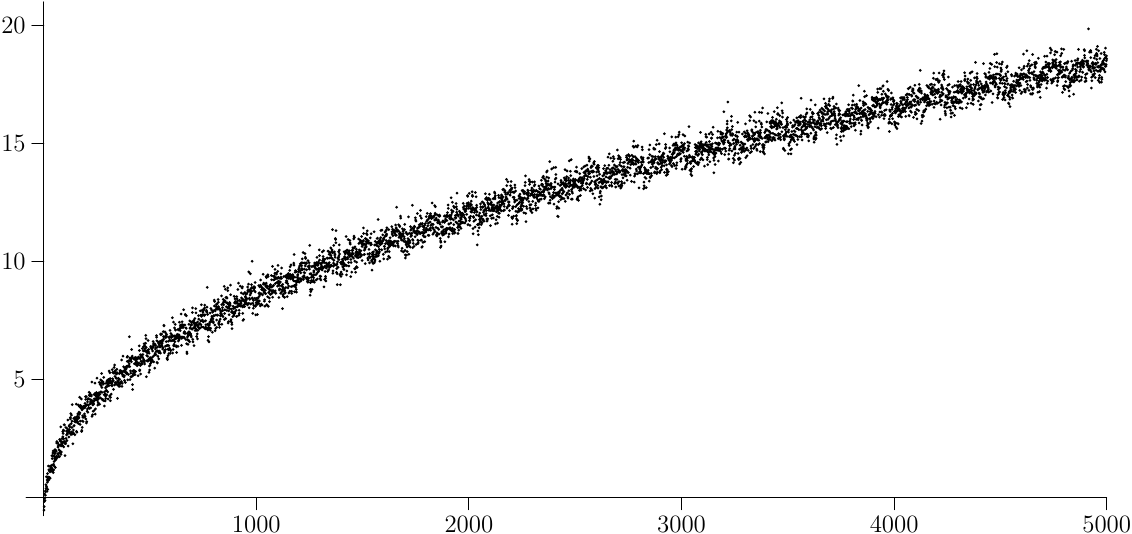}
\caption{Error $\frac{T}{4\pi}\log\frac{T}{2\pi}-\frac{T}{4\pi}-N(T)$ in terms of $n=N(T)$ with $T=\gamma_n$.}
\label{Siegel-rest}
\end{center}
\end{figure}

\section{Number of zeros}\label{S:5}
Let $N(t)$ be the number of zeros of $\Rzeta(s)$ with $0<\Im(\rho)\le t$. Since the zeros in our file are not strictly ordered by increasing imaginary parts, it is not true that if $\rho_n=\beta_n+i\gamma_n$ we have $N(\gamma_n)=n$. But  the difference is only one or two unities for the zeros we have computed. Therefore, we may forget about this difference in the following considerations.

Siegel \cite{Siegel} shows that 
\[N(T)=\frac{T}{4\pi}\log\frac{T}{4\pi}-\frac{T}{4\pi}+\orden(T).\]
We plot the points $(n,\frac{\gamma_n}{4\pi}\log\frac{\gamma_n}{4\pi}-\frac{\gamma_n}{4\pi}-n)$ for $1\le n\le 5000$
we see that in fact there is a good agreement because the error for $N(t)=1000$ is less than $20$. The regularity of the figure makes us think of the possibility of the second term in the approximation.
Experimentally, I find good agreement with the conjecture
\begin{equation}
N(T)=\frac{T}{4\pi}\log\frac{T}{2\pi}-\frac{T}{4\pi}-\frac12\sqrt{\frac{T}{2\pi}}+\frac32+O(\log T).
\end{equation}
Representing the points 
\[\Bigl(n,\frac{\gamma_n}{4\pi}\log\frac{\gamma_n}{2\pi}-\frac{\gamma_n}{4\pi}-\frac12\sqrt{\frac{\gamma_n}{2\pi}}+\frac32-n\Bigr)\]
we obtain the plot
\begin{figure}[H]
\begin{center}
\includegraphics[width=0.6\hsize]{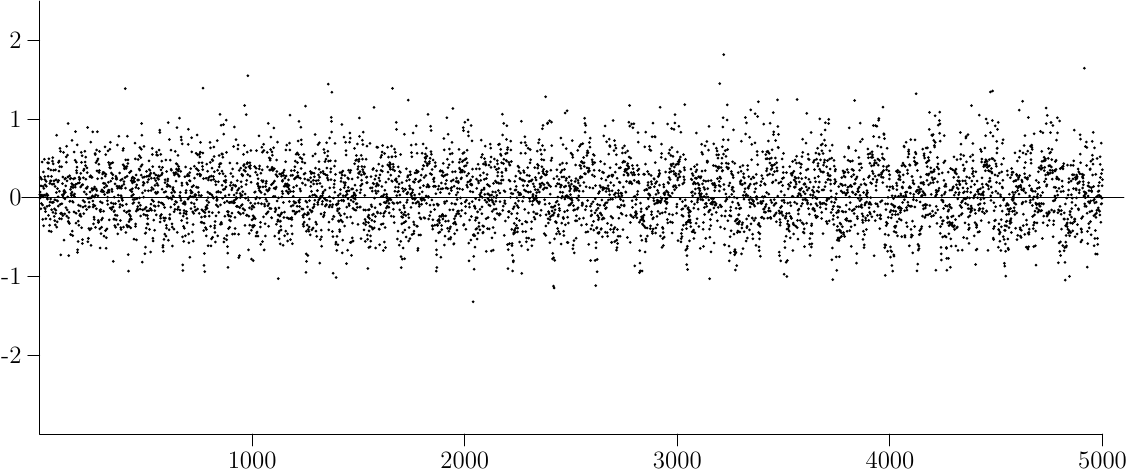}
\caption{Computed error for $n\le 5000$. }
\label{numberzerosD}
\end{center}
\end{figure}

For any real number $\sigma$, denote by $N(\sigma, T)$ the number of zeros $\rho=\beta+i\gamma$ of $\Rzeta(s)$ with $\beta<\sigma$ and 
$0<\gamma<T$. For the set of zeros that I have computed and any $\sigma$ we consider the function 
\[\sum_{\beta_n<\sigma}\Bigl\{
A\Bigl(\frac{\gamma}{4\pi}\log\frac{\gamma}{2\pi}-\frac{\gamma}{4\pi}\Bigr)+B\sqrt{\frac{\gamma}{2\pi}}+C-n\Bigr\}^2.\]
Compute $A$, $B$, and $C$ so that this quantity is minimized. 
We obtain the following values:
\[\begin{array}{ccccccc}
\sigma & A & B & C & m & n & \mu\\
1 &   1.00000001605 & -0.500005726 & 1.43690417259 & 34963.92& 162215 & 0.2155\\
1/2 & 0.99999998564 & -0.499934194 & 1.44317835999 & 22777.01 & 103674 & 0.2196\\
0   & 0.99999950736 & -0.499053472 & 1.52818665442 & 4462.05 & 22983 & 0.1941\\
-1  & 0.99999956381 & -0.499312171 & 1.60544788526 & 1421.67 & 8565 & 0.1659\\
\end{array}\]
$m$ is the min value attained, $n$ is the number of zeros that satisfy the inequality $\beta_n<\sigma$ and $\mu=m/n$ is the mean deviation.

For all computed values of $\rho_n$ we have the bounds for the difference
\[-1.99682<\frac{\gamma_n}{4\pi}\log\frac{\gamma_n}{2\pi}-\frac{\gamma_n}{4\pi}-\frac12\sqrt{\frac{\gamma_n}{2\pi}}+\frac32-n<2.55738\]
Later \cite{A185} we have proved that adding the term  $\frac12\sqrt{\frac{T}{2\pi}}$, the error is at most $\Orden(T^{2/5}\log^2T )$.

\section{Extreme zeros}\label{S:6}
\begin{figure}[H]
\begin{center}
\includegraphics[width=0.8\hsize]{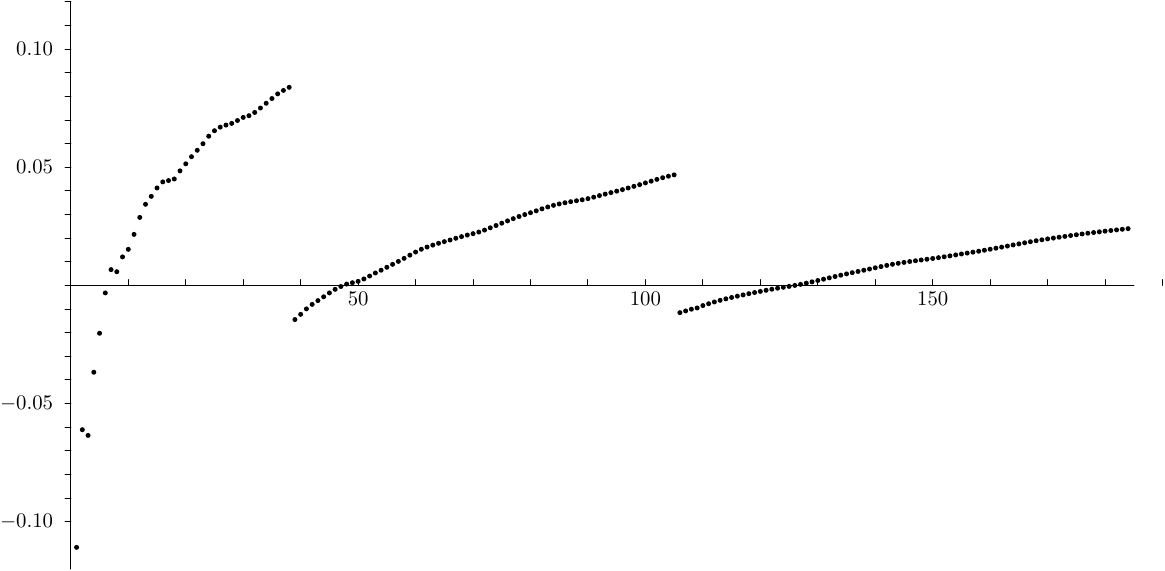}
\caption{Differences $\sqrt{\frac{\gamma_{n_k}}{2\pi}}-(k+1)$}
\label{F:dif}
\end{center}
\end{figure}

The right limit of the zeros with $\rho=\beta+i\gamma$ with $\gamma>0$ is a difficult problem. Siegel stated that for any $\varepsilon>0$ there is some constant $C$ such that $-C\gamma^{\varepsilon}<\beta$, and proved this for $\varepsilon=3/7$. I have been unable to extend his proof beyond  $-C\gamma^{2/5}\log^{4/5}\gamma<\beta$ (\cite{A98}). 
It is natural to consider those zeros $\rho_n$ of our list that are
records for $-\beta_n$. That is such that $-\beta_n>-\beta_m$ for all $1\le m<n$. The corresponding heights $\gamma_n$ are very regular. 
If $\rho_{n_k}$ are the zeros where these records occur, we have 
approximately 
\[\sqrt{\frac{\gamma_{n_k}}{2\pi}}\approx k+1.\]
Between the computed zeros there are $184$ records; in all these cases $k+1$ is the closest integer to $\sqrt{\gamma_k/2\pi}$. 
The differences  also have a remarkable behaviour; see Figure \ref{F:dif}
The jumps occur for $k=38$ where $\sqrt{\frac{\gamma_{n_k}}{2\pi}}-(k+1)=0.0838573$ jumps in the next extreme zero to $-0.0144479$. This corresponds to the zeros numbers 4536 and 4818, with height 9112 and 9597. 

The second jump occurs for  $k=105$ where 
$\sqrt{\frac{\gamma_{n_k}}{2\pi}}-(k+1)=0.0468023$ and jumps to $-0.0115383$. Later, we will  give a possible explanation. In this case, the two extreme zeros are 45792 and 46775 with height 69333 and 70660.

The above means that $\gamma_{n_k}\approx2\pi(k+1)^2$. As shown in \cite{A86} the point $t=2\pi(k+1)^2$ is where a new term $(k+1)^{-s}$  is added to the main approximation of $\Rzeta(s)$. We will see that $\Rzeta(s)$ has a cyclic behaviour associated with these points. 

Our data covers 184 records for the zeros number
\[1,5, 13, 26, 45, 69, 99, 135, 178, 227, 283, 346,\dots, 157671, 159584,161510.\]

The next figure contains the plot of the points $(\gamma,-\beta)$.
Most of the zeros are in the strip  $-1\le\beta\le1$. All  zeros satisfies $\beta\le 1$. We conjecture that this is true for all zeros. In \cite{A173} we show that this is true for zeros with $\gamma\ge t_0$ for a very great $t_0$. 

We see that the  zeros with $-\beta$ large are organized in lines. The jumps in Figure \ref{F:dif} correspond to the points where the line giving the extreme zeros  changes into another. This happens twice  on our sample. And we see that this happens at the height where the jumps take place.

\begin{figure}[H]
\begin{center}
\includegraphics[width=\hsize]{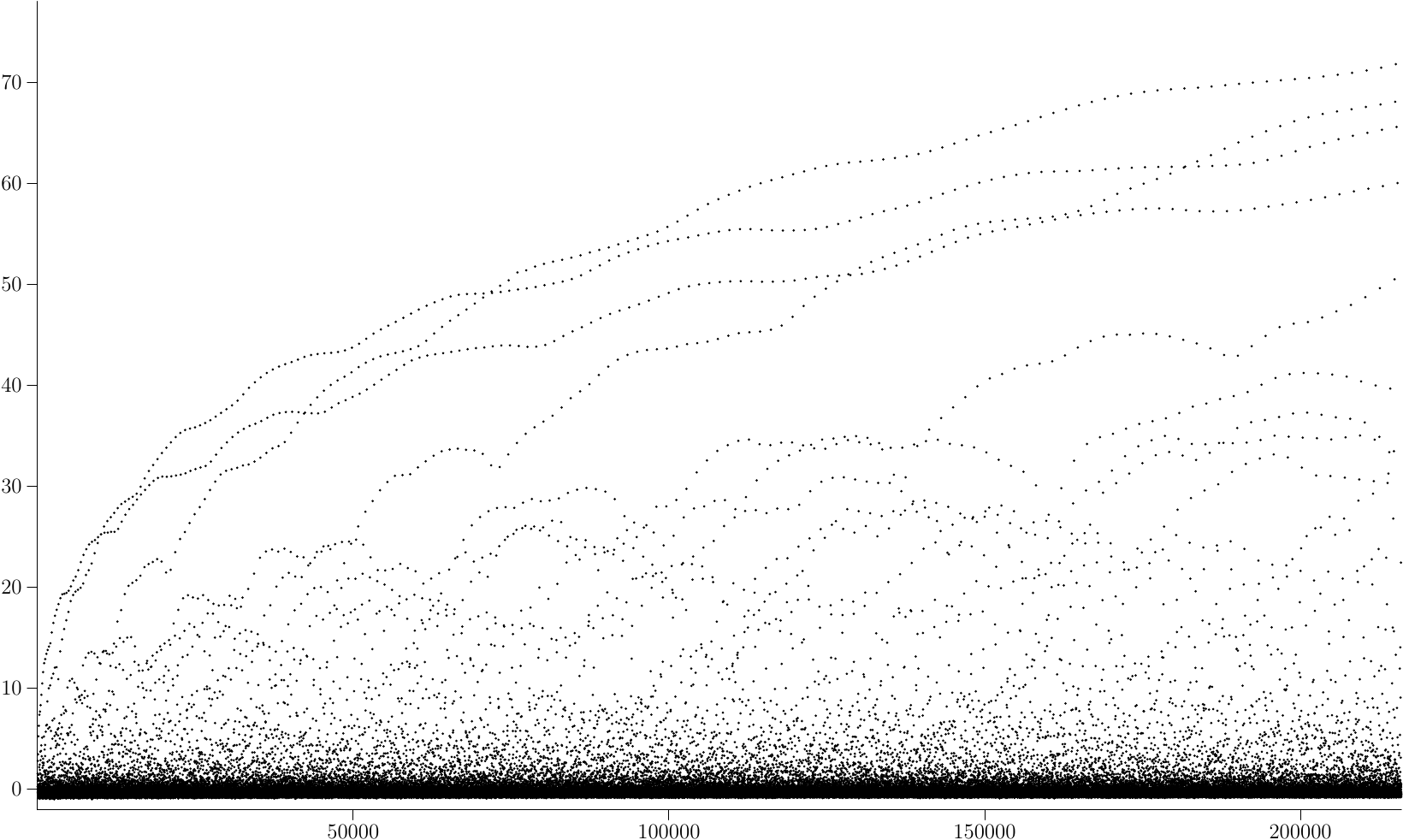}
\caption{Points $(\gamma,-\beta)$}
\label{F:plotzeros2}
\end{center}
\end{figure}

\subsection{Left limit of extremal zeros}
We have found that the sequence of extremal zeros $\rho_{n_k}=\beta_{n_k}+i\gamma_{n_k}$ with $-\beta_{n_k}>-\beta_n$ for all $n<n_k$, have very regular heights.  We are especially interested for the $\beta_{n_k}$. 

I find that $-\beta_{n_k}\approx a(k+1)^{2/3}$. Taking into account the relationship between $k$ and $\gamma_{n_k}$ this means that 

\begin{equation}
-\beta_{n_k}\approx 2.2430\Bigl(\frac{\gamma_{n_k}}{2\pi}\Bigr)^{1/3}.
\end{equation}
But there is little data to get this relation. The errors oscillate between $-1$ and $1.5$ but not in the form of random noise, but rather form certain curves.

\section{Cyclic structure of the zeros}\label{S:7}

Siegel \cite{Siegel}*{p.~308} proved that 
\begin{equation}
h(t):=-\sum_{0<\gamma\le t}\beta=\frac{t}{4\pi}\log2+\orden(t),\qquad t\to+\infty.
\end{equation}

Some structure in the zeros of $\Rzeta(s)$ is not shown in Figure \ref{F:plotzeros2}. There is a cyclic structure that repeats in each section 
$[2\pi n^2, 2\pi(n+1)^2)$. To show this, we show two plots.
The first shows the points $(\gamma,-\beta)$ for the zeros $\rho=\beta+i\gamma$ with $2\pi60^2\le\gamma<2\pi65^2$. Therefore, we plot $5$ cycles 
\begin{figure}[H]
\begin{center}
\includegraphics[width=\hsize]{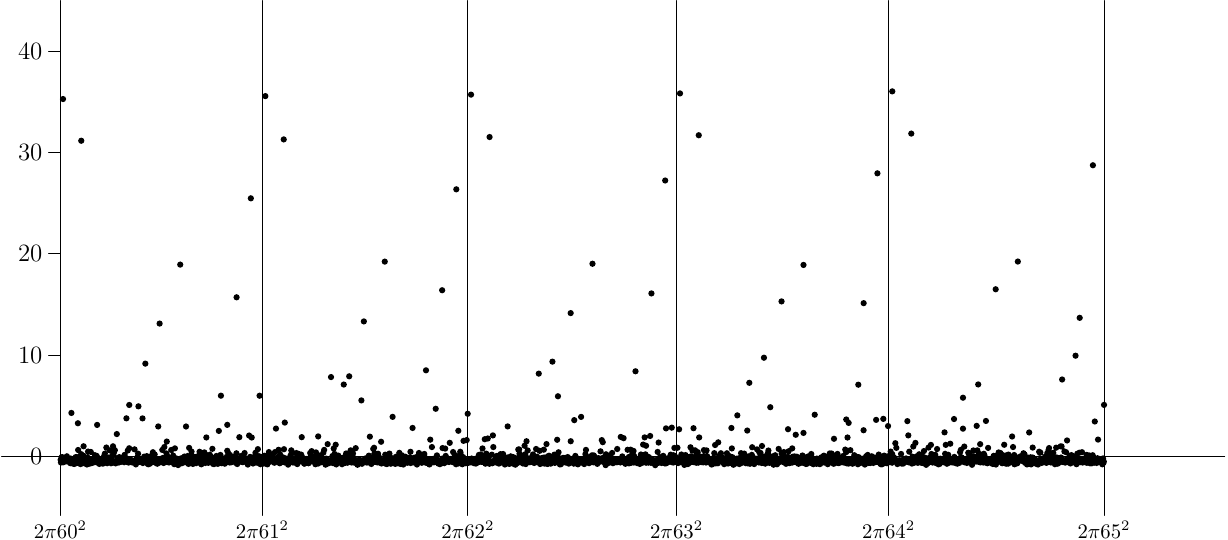}
\caption{Zeros with $2\pi60^2\le\gamma<2\pi65^2$}
\label{F:cycles1}
\end{center}
\end{figure}
It appears as if we can recognize some of the zeros in the different cycles. For example, after the zero with greatest $-\beta$ we see another one a little lesser with a greater $\gamma$. It seems that these structures evolve. We may compare with the next figure for the zeros with $2\pi179^2\le\gamma<2\pi185^2$
\begin{figure}[H]
\begin{center}
\includegraphics[width=\hsize]{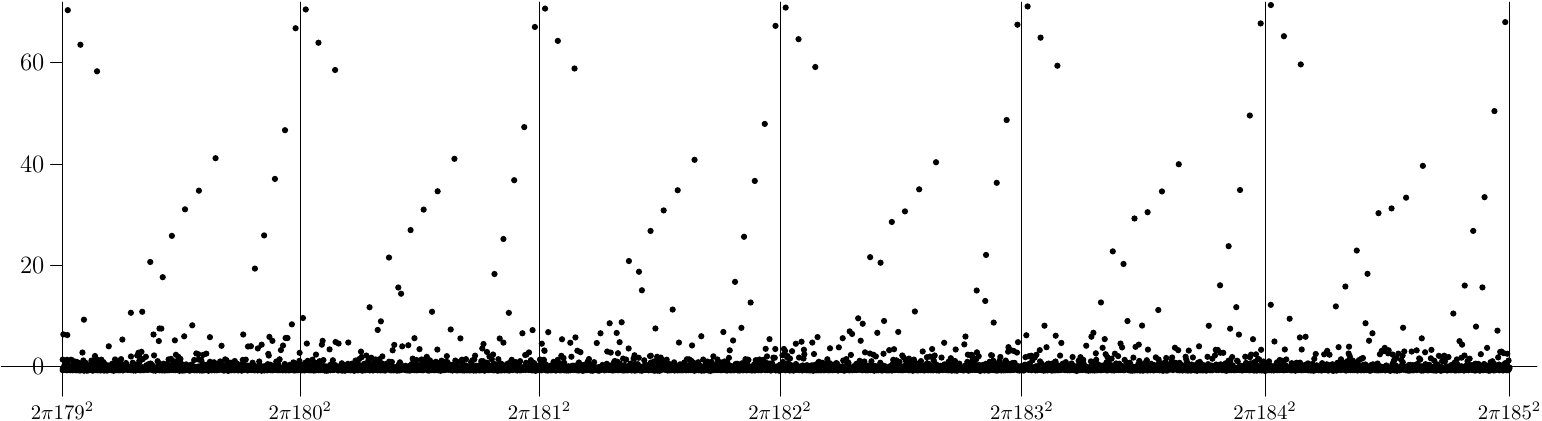}
\caption{Zeros with $2\pi179^2\le\gamma<2\pi185^2$}
\label{F:cycles2}
\end{center}
\end{figure}

\vspace{12pt}

\section{Siegel sum}\label{S:8}
One of Siegel's Theorems in \cite{Siegel} is
\[h(T):=-\sum_{\gamma\le T}\beta=\frac{T}{4\pi}\log2+\orden(T).\]
We see that $h(x)$ has a  cyclic behaviour,  related to the cyclic behaviour of the zeros. We also see that Siegel's approximation is very good, so we plot the difference. For the first 100000 zeros, we get 

\begin{figure}[H]
\begin{center}
\includegraphics[width=\hsize]{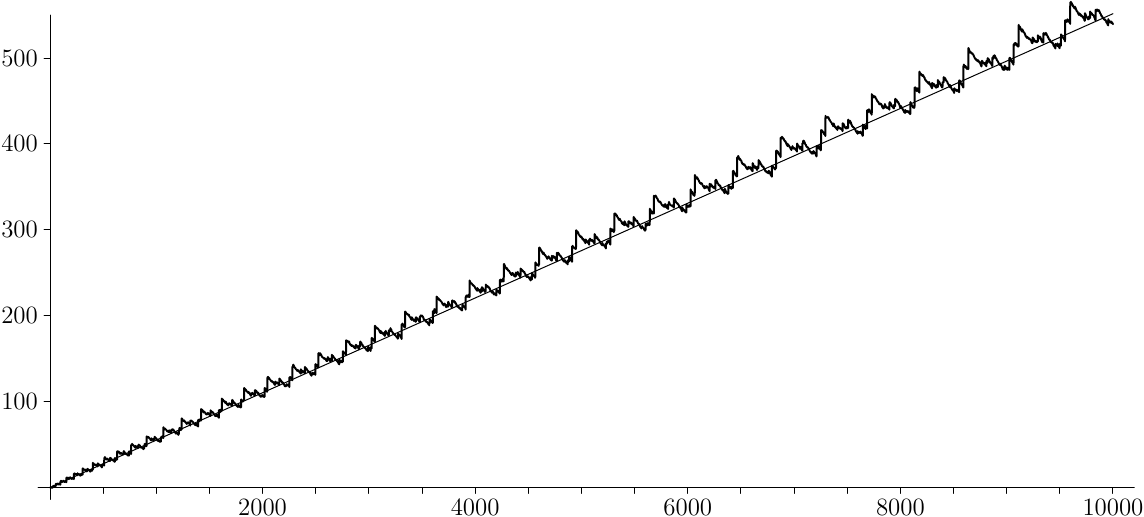}
\caption{Plot of $h(T)$ and its approximation $\frac{T}{4\pi}\log2$}
\label{Siegelfunction1}
\end{center}
\end{figure}
\begin{figure}[H]
\begin{center}
\includegraphics[width=\hsize]{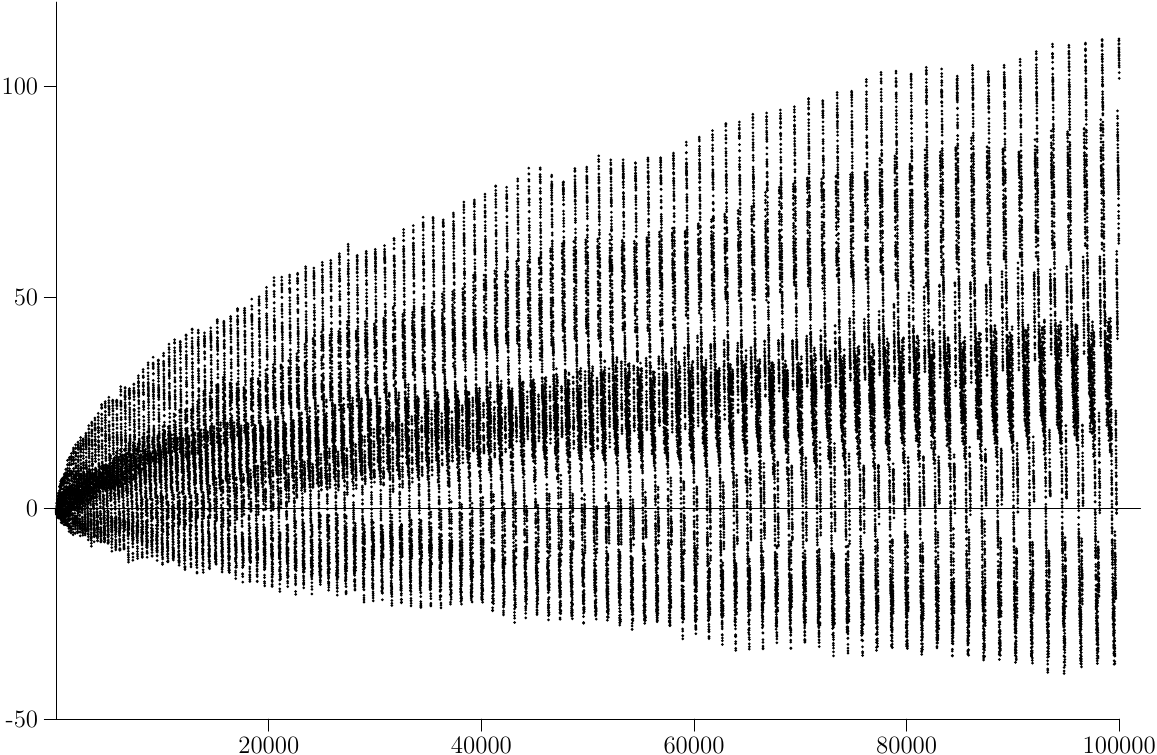}
\caption{Plot of $h(t)-\frac{t}{4\pi}\log2$}
\label{diferenciah}
\end{center}
\end{figure}
To see the cyclic behavior, we plot this difference in three contiguous cycles for 
$2\pi 100^2<t<2\pi 103^2$
\begin{figure}[H]
\begin{center}
\includegraphics[width=0.8\hsize]{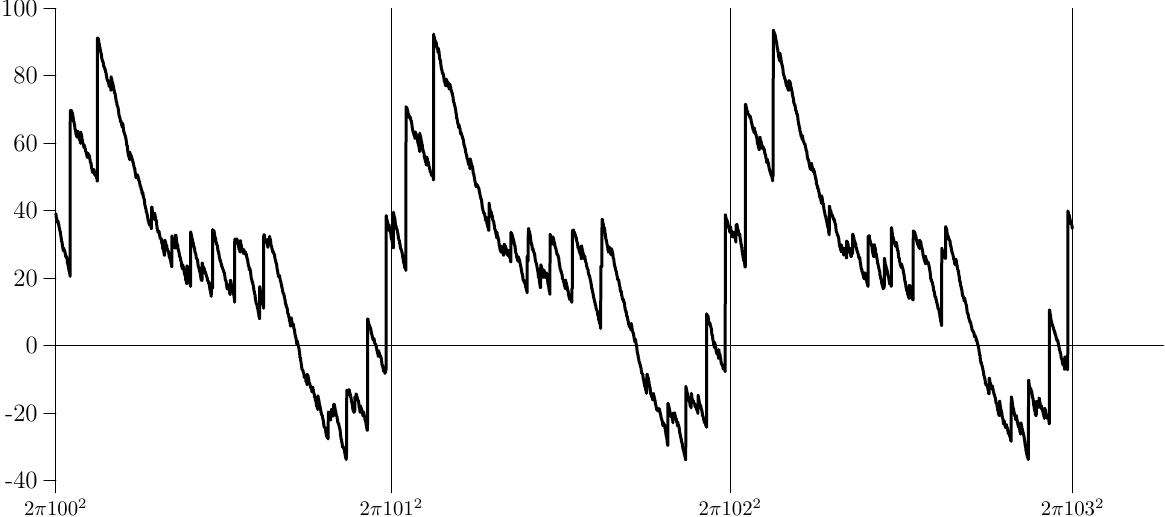}
\caption{Plot of $h(t)-\frac{t}{4\pi}\log2$ in $2\pi 100^2<t<2\pi 103^2$}
\label{difhD}
\end{center}
\end{figure}
To compare, we show another section.
\begin{figure}[H]
\begin{center}
\includegraphics[width=0.8\hsize]{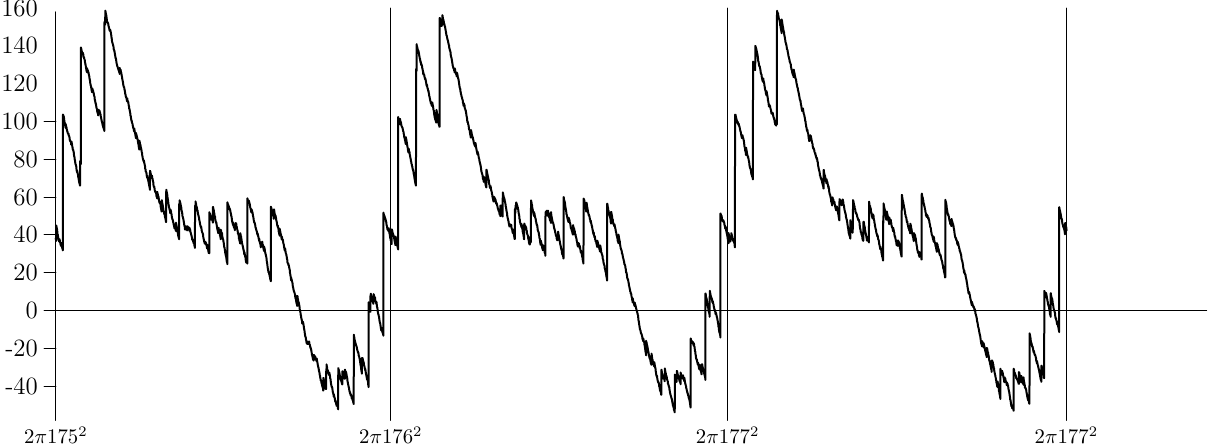}
\caption{Plot of $h(t)-\frac{t}{4\pi}\log2$ in $2\pi 175^2<t<2\pi 178^2$}
\label{difhU}
\end{center}
\end{figure}

In a previous version of \cite{A102} we proved 
\begin{lemma}
Let $t_0$ be a fixed real number. For $T\to+\infty$ we have
\[\Im\int_{1+t_0}^{1+iT}\log\frac{1-e^{2\pi i\eta}}{1+e^{4\pi i\eta}}\,ds=-\sqrt{T/2\pi}P(\sqrt{T/2\pi})+\Orden(1).\]
where
\[P(x)=\sum_{n=1}^\infty\frac{2\sin2\pi n x}{n^2}-\sum_{n=1}^\infty(-1)^n\frac{\sin(4\pi n x)}{n^2}.\]
\end{lemma}
This term will be in the equation for $-\sum_{0<\gamma\le T}\beta$ except that greater error terms eliminate it from the final result. It is because of this that this computation has been eliminated in the final form of \cite{A102}. But it appears to explain Figures \ref{difhD} and \ref{difhU}.

Hence, we have the conjecture
\begin{equation}\label{E:conj}
h(t)=\frac{t}{4\pi}\log 2-\Bigl(\frac14+\frac{P(\sqrt{t/2\pi})}{2\pi}\Bigr)\sqrt{\frac{t}{2\pi}}+\Orden(t^{2/5}\log t).
\end{equation}
Although the difference has jumps of size $\beta$ in each zero $\rho=\beta+i\gamma$, that presumably are not bounded. We put $\Orden(t^{2/5}\log t)$ which is the bound of the $\beta$ that we know. 

\begin{figure}[H]
\begin{center}
\includegraphics[width=0.8\hsize]{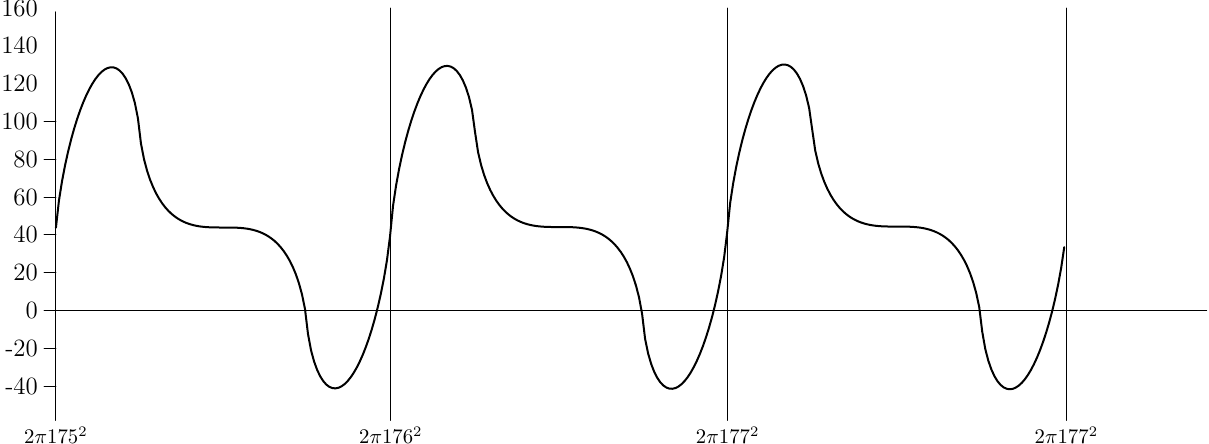}
\caption{Plot of $\Bigl(\frac14+\frac{P(\sqrt{t/2\pi})}{2\pi}\Bigr)\sqrt{\frac{t}{2\pi}}$}
\label{extraterm}
\end{center}
\end{figure}

Compare Figure \ref{difhU} with Figure \ref{extraterm}
The corresponding error is represented here
\begin{figure}[H]
\begin{center}
\includegraphics[width=0.8\hsize]{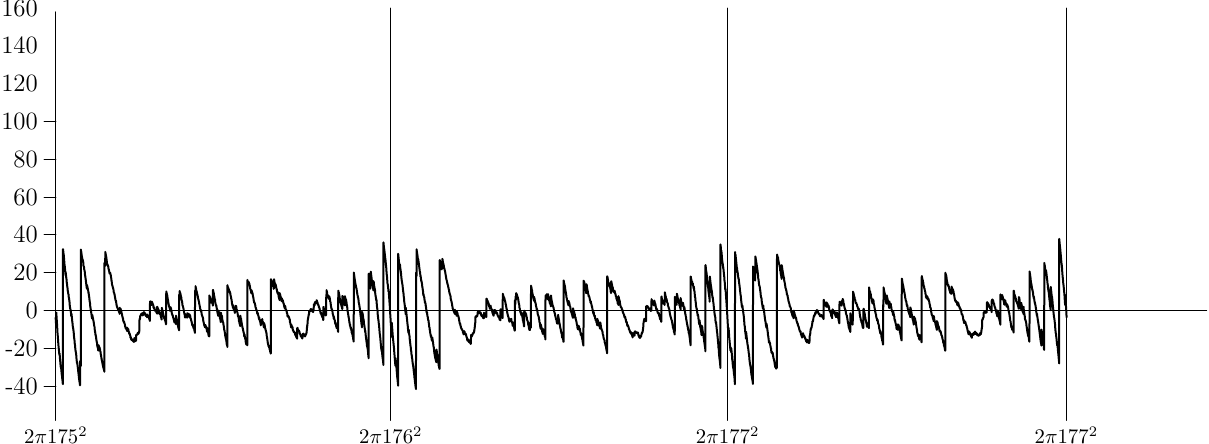}
\caption{Error $h(t)-\frac{t}{4\pi}\log 2-\Bigl(\frac14+\frac{P(\sqrt{t/2\pi})}{2\pi}\Bigr)\sqrt{\frac{t}{2\pi}}$}
\label{extratermerror}
\end{center}
\end{figure}

\section{Horizontal distribution of zeros}\label{S:9}
\begin{figure}[H]
\begin{center}
\includegraphics[width=0.8\hsize]{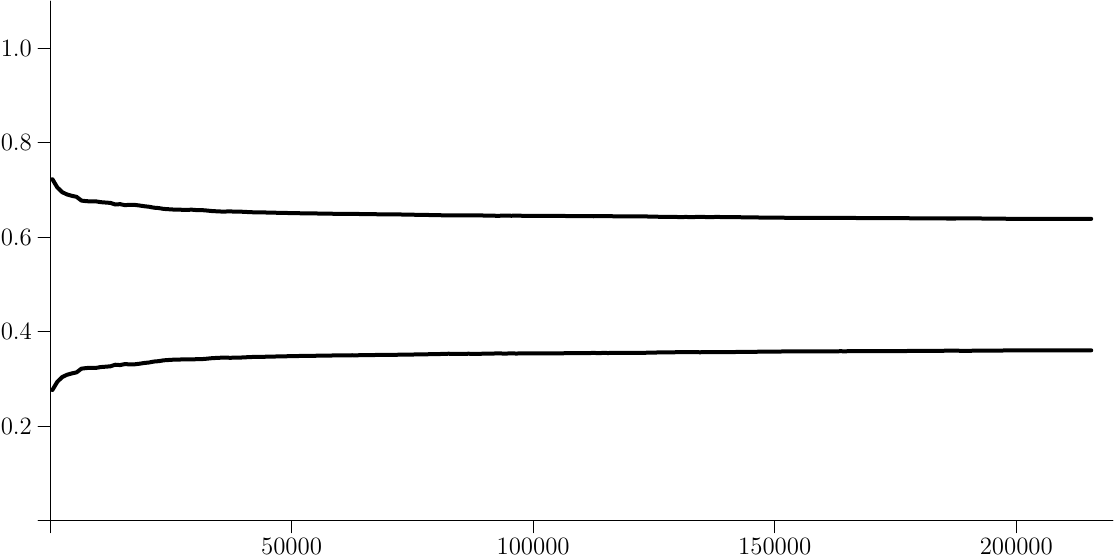}
\caption{Evolution of the density of zeros}
\label{densfig}
\end{center}
\end{figure}

\begin{figure}[H]
\begin{center}
\includegraphics[width=\hsize]{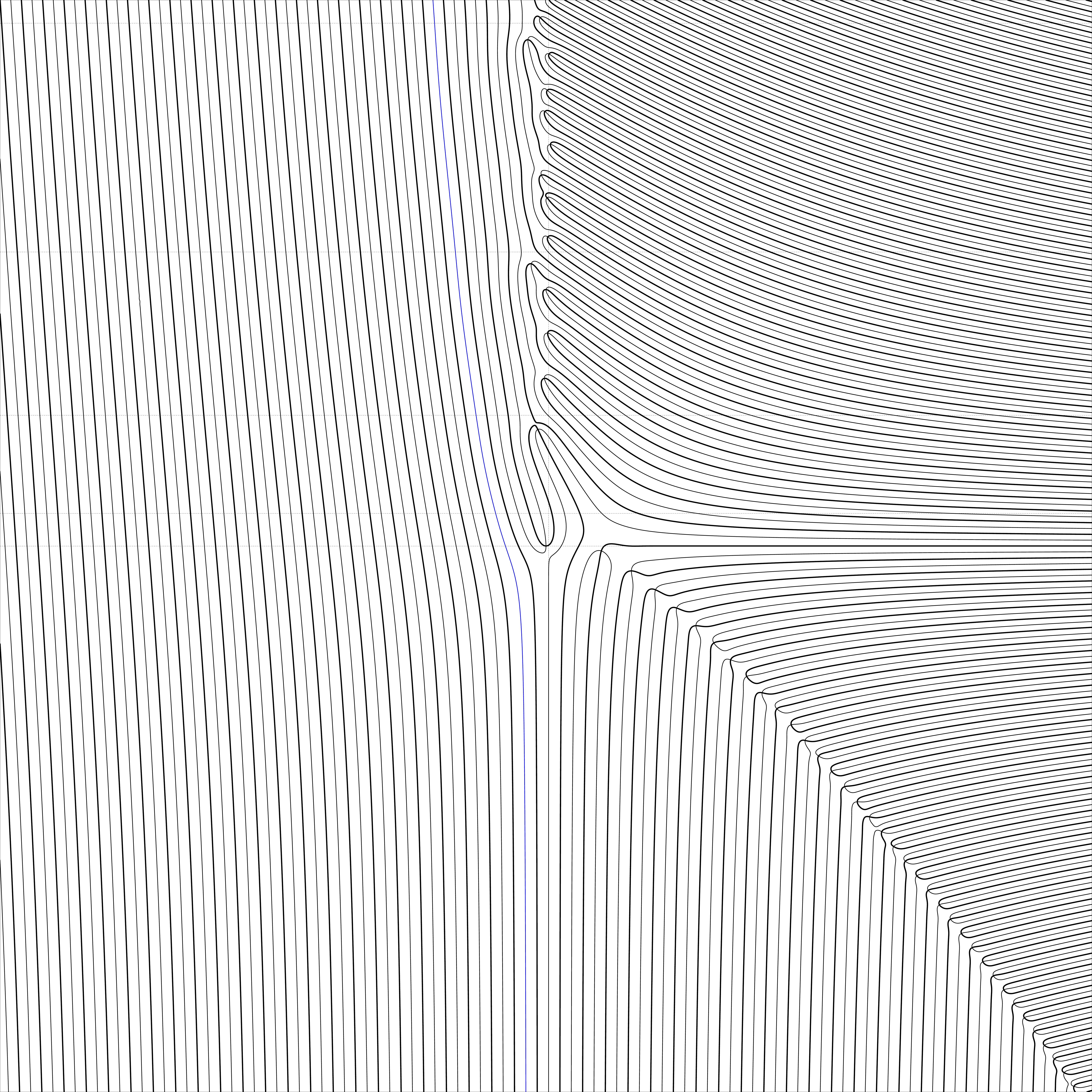}
\caption{x-ray of $\pi^{-s/2}\Gamma(s/2)\Rzeta(s)$ in $(-105,105)^2$, horizontal lines at $2\pi n^2$.}
\label{Rzetapeinada}
\end{center}
\end{figure}

As we see in the x-ray of $\pi^{-s/2}\Gamma(s/2)\Rzeta(s)$ (see Figure \ref{Rzetapeinada}) all the imaginary lines containing the zeros $\rho=\beta+i\gamma$ of $\Rzeta(s)$ with $\gamma>0$ goes to the right. Therefore, each zero with $\beta<\frac12$ is associated with two zeros of $\zeta(s)$ at the points where the imaginary line cuts the critical line. 
Therefore, if a proportion $\delta$ of the $\frac{T}{4\pi}\log\frac{T}{2\pi}-\frac{T}{4\pi}$ zeros of $\Rzeta(s)$ with $0<\gamma\le T$ satisfies $\beta<1/2$, then a proportion $\delta$ of zeros of $\zeta(s)$ will be on the critical line. 
Therefore, it is interesting to compute the proportion of zeros of $\Rzeta(s)$ with $\beta<\frac12$. From the 162215 computed zeros, 103674 satisfies $\beta<1/2$ and 58541 $\beta>1/2$. So, the total densities are 
\[\frac{103674}{162215}=0.639115,\quad \frac{58541}{162215}=0.360885\]
The first of these densities appear to be decreasing, and the other to be increasing, but only very slightly. Figure \ref{densfig} represents the graphic of 
\[\delta(t)=\frac{|\{\rho=\beta+i\gamma\colon 0<\gamma<t, \beta<1/2\}|}
{|\{\rho=\beta+i\gamma\colon 0<\gamma<t\}|}\]
and below the graph of $1-\delta(t)$, the density of zeros with $\beta\ge1/2$.  

Fixing a natural number $N$, we separate the real line in intervals $[k/N,(k+1)/N)$ and compute the number of zeros $a(k)$ with $\frac{k}{N}\le \beta<\frac{k+1}{N}$. Then we plot the function that in this interval takes the value 
$d(k)=\frac{a(k)N}{Z}$, where $Z$ is the total number of zeros in our list $Z=162265$. We get in this way, with $N=26$ the histogram
This histogram takes all the zeros $\rho=\beta+i\gamma$ with $0<\gamma\le215\,946.3$. We can determine in the same way the histograms $H(T,M)$ with all zeros satisfying $0<\gamma\le T$ and with divisions of length $1/M$.  The results in \cite{A174} implies that  for each fixed $N$ and with $T\to+\infty$ the heights for all intervals $[k/M, (k+1)/M)$ with $k/M>1/2$ tends to $0$. We do not know what happens to the intervals to the right of $1/2$. It will be interesting to compute a set of higher zeros to see if there is any trend.

In an appendix, we put a  table with the zero count in each section $2\pi (n-1)^2\le\gamma<2\pi n^2$ classified according with $\beta$ in the intervals 
$[1/2,1)$, $(-\infty,1/2)$, $[0,1/2)$ and $(-\infty,0)$.

These proportions, if proved for $T\to+\infty$,  will be a great advance over the results obtained with the Levinson method.

\newpage

\setlength{\tabcolsep}{3pt}

\begin{tabular}{|c|c|c|c|c|}
\hline 
$n$ & $\beta\in[\frac12, 1)$ & $\beta<\frac12$ & $\beta\in[0,\frac12)$ & $\beta<0$\\ \hline
1 & 0 & 0 & 0 & 0\\ 
2 & 0 & 1 & 0 & 1\\ 
3 & 0 & 4 & 2 & 2\\ 
4 & 2 & 6 & 3 & 3\\ 
5 & 4 & 9 & 4 & 5\\ 
6 & 7 & 12 & 5 & 7\\ 
7 & 8 & 16 & 8 & 8\\ 
8 & 5 & 23 & 15 & 8\\ 
9 & 11 & 25 & 16 & 9\\ 
10 & 14 & 29 & 18 & 11\\ 
11 & 13 & 36 & 23 & 13\\ 
12 & 18 & 37 & 19 & 18\\ 
13 & 19 & 44 & 27 & 17\\ 
14 & 23 & 47 & 28 & 19\\ 
15 & 22 & 55 & 36 & 19\\ 
16 & 25 & 59 & 37 & 22\\ 
17 & 31 & 61 & 42 & 19\\ 
18 & 29 & 71 & 45 & 26\\ 
19 & 30 & 78 & 52 & 26\\ 
20 & 42 & 73 & 47 & 26\\ 
21 & 46 & 77 & 53 & 24\\ 
22 & 37 & 95 & 65 & 30\\ 
23 & 43 & 97 & 68 & 29\\ 
24 & 47 & 101 & 63 & 38\\ 
25 & 48 & 108 & 71 & 37\\ 
26 & 48 & 116 & 84 & 32\\ 
27 & 63 & 111 & 71 & 40\\ 
28 & 54 & 127 & 92 & 35\\ 
29 & 62 & 129 & 87 & 42\\ 
30 & 68 & 131 & 89 & 42\\ 
31 & 68 & 140 & 96 & 44\\ 
32 & 84 & 133 & 92 & 41\\ 
33 & 74 & 152 & 107 & 45\\ 
34 & 83 & 152 & 103 & 49\\ 
35 & 80 & 164 & 115 & 49\\ 
36 & 74 & 179 & 123 & 56\\ 
37 & 90 & 172 & 122 & 50\\ 
38 & 88 & 183 & 129 & 54\\ 
39 & 94 & 186 & 126 & 60\\ 
40 & 98 & 192 & 141 & 51\\ 
41 & 100 & 200 & 142 & 58\\ 
42 & 101 & 208 & 151 & 57\\ 
43 & 108 & 210 & 155 & 55\\ 
44 & 107 & 221 & 166 & 55\\ 
45 & 117 & 220 & 158 & 62\\ 
46 & 129 & 218 & 154 & 64\\ \hline
\end{tabular}
\quad
\begin{tabular}{|c|c|c|c|c|}
\hline 
$n$ & $\beta\in[\frac12, 1)$ & $\beta<\frac12$ & $\beta\in[0,\frac12)$ & $\beta<0$\\ \hline
47 & 124 & 232 & 172 & 60\\ 
48 & 116 & 250 & 189 & 61\\ 
49 & 137 & 238 & 174 & 64\\ 
50 & 126 & 260 & 199 & 61\\ 
51 & 135 & 261 & 198 & 63\\ 
52 & 130 & 276 & 201 & 75\\ 
53 & 145 & 270 & 203 & 67\\ 
54 & 144 & 281 & 210 & 71\\ 
55 & 157 & 279 & 208 & 71\\ 
56 & 155 & 290 & 213 & 77\\ 
57 & 163 & 292 & 211 & 81\\ 
58 & 181 & 285 & 212 & 73\\ 
59 & 160 & 316 & 227 & 89\\ 
60 & 182 & 303 & 223 & 80\\ 
61 & 178 & 318 & 239 & 79\\ 
62 & 185 & 322 & 239 & 83\\ 
63 & 183 & 332 & 251 & 81\\ 
64 & 185 & 342 & 260 & 82\\ 
65 & 189 & 349 & 253 & 96\\ 
66 & 191 & 356 & 274 & 82\\ 
67 & 185 & 372 & 283 & 89\\ 
68 & 201 & 368 & 275 & 93\\ 
69 & 203 & 375 & 289 & 86\\ 
70 & 205 & 384 & 287 & 97\\ 
71 & 204 & 395 & 299 & 96\\ 
72 & 236 & 375 & 266 & 109\\ 
73 & 229 & 392 & 298 & 94\\ 
74 & 231 & 400 & 295 & 105\\ 
75 & 230 & 412 & 310 & 102\\ 
76 & 221 & 431 & 345 & 86\\ 
77 & 228 & 436 & 331 & 105\\ 
78 & 238 & 435 & 329 & 106\\ 
79 & 255 & 430 & 326 & 104\\ 
80 & 236 & 459 & 350 & 109\\ 
81 & 259 & 447 & 328 & 119\\ 
82 & 258 & 458 & 345 & 113\\ 
83 & 256 & 472 & 355 & 117\\ 
84 & 260 & 479 & 369 & 110\\ 
85 & 273 & 476 & 359 & 117\\ 
86 & 268 & 493 & 380 & 113\\ 
87 & 286 & 484 & 371 & 113\\ 
88 & 268 & 515 & 391 & 124\\ 
89 & 284 & 508 & 386 & 122\\ 
90 & 293 & 512 & 387 & 125\\ 
91 & 292 & 522 & 406 & 116\\ 
92 & 288 & 539 & 418 & 121\\ \hline
\end{tabular}
\label{table1}

\newpage

\renewcommand{\arraystretch}{0.9}

\begin{tabular}{|c|c|c|c|c|}
\hline 
$n$ & $\beta\in[\frac12, 1)$ & $\beta<\frac12$ & $\beta\in[0,\frac12)$ & $\beta<0$\\ \hline
93 & 297 & 539 & 417 & 122\\ 
94 & 305 & 544 & 419 & 125\\ 
95 & 301 & 558 & 435 & 123\\ 
96 & 320 & 550 & 424 & 126\\ 
97 & 320 & 562 & 432 & 130\\ 
98 & 315 & 577 & 451 & 126\\ 
99 & 319 & 585 & 452 & 133\\ 
100 & 321 & 594 & 456 & 138\\ 
101 & 335 & 591 & 460 & 131\\ 
102 & 347 & 591 & 454 & 137\\ 
103 & 329 & 618 & 481 & 137\\ 
104 & 338 & 623 & 485 & 138\\ 
105 & 358 & 613 & 485 & 128\\ 
106 & 346 & 636 & 496 & 140\\ 
107 & 356 & 638 & 496 & 142\\ 
108 & 361 & 644 & 502 & 142\\ 
109 & 368 & 649 & 489 & 160\\ 
110 & 383 & 644 & 496 & 148\\ 
111 & 369 & 671 & 533 & 138\\ 
112 & 389 & 663 & 510 & 153\\ 
113 & 384 & 678 & 533 & 145\\ 
114 & 401 & 672 & 524 & 148\\ 
115 & 399 & 686 & 528 & 158\\ 
116 & 379 & 718 & 559 & 159\\ 
117 & 392 & 716 & 567 & 149\\ 
118 & 408 & 711 & 558 & 153\\ 
119 & 399 & 733 & 573 & 160\\ 
120 & 424 & 718 & 563 & 155\\ 
121 & 425 & 729 & 579 & 150\\ 
122 & 414 & 753 & 600 & 153\\ 
123 & 427 & 750 & 592 & 158\\ 
124 & 413 & 776 & 614 & 162\\ 
125 & 438 & 763 & 596 & 167\\ 
126 & 433 & 779 & 614 & 165\\ 
127 & 437 & 788 & 617 & 171\\ 
128 & 438 & 796 & 623 & 173\\ 
129 & 456 & 791 & 618 & 173\\ 
130 & 446 & 814 & 637 & 177\\ 
131 & 462 & 808 & 629 & 179\\ 
132 & 466 & 818 & 636 & 182\\ 
133 & 467 & 826 & 649 & 177\\ 
134 & 472 & 835 & 645 & 190\\ 
135 & 476 & 842 & 659 & 183\\ 
136 & 476 & 854 & 681 & 173\\ 
137 & 499 & 842 & 665 & 177\\ 
138 & 473 & 881 & 696 & 185\\
139 & 481 & 884 & 692 & 192\\   \hline
\end{tabular}
\quad
\begin{tabular}{|c|c|c|c|c|}
\hline 
$n$ & $\beta\in[\frac12, 1)$ & $\beta<\frac12$ & $\beta\in[0,\frac12)$ & $\beta<0$\\ \hline
140 & 513 & 865 & 658 & 207\\ 
141 & 523 & 866 & 690 & 176\\ 
142 & 549 & 853 & 655 & 198\\ 
143 & 506 & 906 & 725 & 181\\ 
144 & 531 & 894 & 704 & 190\\ 
145 & 513 & 923 & 723 & 200\\ 
146 & 517 & 932 & 734 & 198\\ 
147 & 524 & 937 & 747 & 190\\ 
148 & 524 & 948 & 750 & 198\\ 
149 & 536 & 949 & 749 & 200\\ 
150 & 572 & 926 & 743 & 183\\ 
151 & 552 & 955 & 747 & 208\\ 
152 & 567 & 954 & 765 & 189\\ 
153 & 581 & 952 & 743 & 209\\ 
154 & 581 & 964 & 756 & 208\\ 
155 & 575 & 982 & 780 & 202\\ 
156 & 580 & 989 & 793 & 196\\ 
157 & 578 & 1004 & 799 & 205\\ 
158 & 567 & 1025 & 816 & 209\\ 
159 & 590 & 1016 & 807 & 209\\ 
160 & 598 & 1019 & 801 & 218\\ 
161 & 587 & 1043 & 826 & 217\\ 
162 & 588 & 1054 & 843 & 211\\ 
163 & 606 & 1048 & 841 & 207\\ 
164 & 620 & 1046 & 827 & 219\\ 
165 & 611 & 1068 & 839 & 229\\ 
166 & 606 & 1084 & 864 & 220\\ 
167 & 633 & 1070 & 840 & 230\\ 
168 & 635 & 1079 & 866 & 213\\ 
169 & 639 & 1088 & 868 & 220\\ 
170 & 638 & 1103 & 876 & 227\\ 
171 & 643 & 1108 & 878 & 230\\ 
172 & 657 & 1108 & 878 & 230\\ 
173 & 642 & 1133 & 916 & 217\\ 
174 & 636 & 1154 & 921 & 233\\ 
175 & 652 & 1148 & 914 & 234\\ 
176 & 699 & 1115 & 864 & 251\\ 
177 & 677 & 1148 & 923 & 225\\ 
178 & 692 & 1147 & 915 & 232\\ 
179 & 674 & 1176 & 954 & 222\\ 
180 & 689 & 1174 & 931 & 243\\ 
181 & 683 & 1192 & 953 & 239\\ 
182 & 660 & 1227 & 973 & 254\\ 
183 & 703 & 1198 & 958 & 240\\ 
184 & 689 & 1223 & 966 & 257\\ 
185 & 704 & 1221 & 986 & 235\\ 
186 & 272 & 480 & 382 & 98\\ \hline
\end{tabular}

\end{document}